\documentclass[addpoints]{amsart} 
\usepackage[hidelinks]{hyperref} 
\hypersetup{allcolors=black}

\usepackage{gensymb} 
\usepackage{pgf, tikz} 
\usepackage{tikz-cd} 
\usetikzlibrary{angles,arrows.meta,automata,backgrounds,calc,decorations.markings,decorations.pathreplacing,intersections,patterns,positioning,quotes} 
\usetikzlibrary{shapes} 
\usepgflibrary{shapes.geometric}
\usepackage{tkz-euclide} 
\usepackage{etoolbox}
\usepackage{placeins}

\usepackage{mdframed} 
\usepackage{pgfplots}
\usepackage{tikz-3dplot} 
\pgfplotsset{width=10cm,compat=1.9}
\usepackage[export]{adjustbox}
\usepackage{caption}
\usepackage{subcaption}
\usepackage{wrapfig}
\usepackage{graphicx}
\usepackage{mathtools}  
\usepackage{amsthm,thmtools} 
\usepackage{relsize} 
\usepackage{array}   
\newcolumntype{L}{>{$}l<{$}} 

\usepackage{amsmath,systeme} 

\usepackage{fancyhdr} 

\usepackage{fdsymbol} 
\usepackage{dutchcal} 
\usepackage{courier}
\usepackage[scr=rsfso]{mathalfa} 
\input ArtNouvc.fd
\input Starburst.fd

\usepackage[all,cmtip]{xy}

\usepackage{enumerate}
\usepackage{multicol} 
\usepackage{mdwlist}

\usepackage[OT2,T1]{fontenc}
\DeclareSymbolFont{cyrletters}{OT2}{wncyr}{m}{n}
\DeclareMathSymbol{\Sha}{\mathalpha}{cyrletters}{"58}

\theoremstyle{plain}
\newtheorem{Theorem}[subsection]{Theorem}
\newtheorem{Lemma}[subsection]{Lemma}
\newtheorem{Corollary}[subsection]{Corollary}

\newtheorem{Definition}[subsection]{Definition}

\theoremstyle{remark}
\newtheorem{Example}[subsubsection]{\bf Example}

\newtheorem{Remark}[subsubsection]{\bf Remark}
\newtheorem{Remarks}[subsubsection]{\bf Remarks}
\numberwithin{equation}{subsection}

\newcommand{\bs}[1]{{\boldsymbol{#1}}}
\newcommand{\scr}[1]{{\mathscr{#1}}}

\newcommand{\mbb}[1]{{\mathbb{#1}}}

\renewcommand{\a}{{\bs a}}
\renewcommand{\b}{{\bs b}}
\renewcommand{\c}{{\bs c}}

\renewcommand{\i}{{\boldsymbol i}}

\newcommand\x{{\bold x}}

\newcommand{\C}{{\mathbb C}}

\newcommand{\FF}{{\scr{F}}}

\renewcommand{\P}{{\mathbb P}}

\newcommand{\R}{{\mathbb R}}

\newcommand{\T}{{\mathsf{T}}}

\newcommand{\Z}{{\mathbb Z}}

\newcommand{\br}[1]{{\langle{#1}\rangle}}

\newcommand{\df}{{\,\overset{\text{\rm df}}{=}\,}}

\newcommand{\ep}{{\varepsilon}}

\newcommand{\id}{{\text{\rm id}}}

\newcommand{\lm}{\longmapsto}
\newcommand{\lr}{\longrightarrow}

\newcommand{\mywedge}[1]{{\bigwedge}^{\kern-1pt{#1}}}

\renewcommand{\pmod}[1]{{\text{\rm (mod ${#1}$)}}}

\newcommand{\Arc}{{\sf{Arc}}}

\newcommand{\Arg}{{\sf{Arg}}}

\newcommand{\Hom}{{\text{\rm Hom}}}

\newcommand{\Set}{\mathsf{Set}}

\newcommand{\Var}[1]{\mathsf{#1}\text{-}\mathsf{Var}}

\usepackage{placeins}

\renewcommand{\T}{{\mbb T}}

\author[E. Brussel]{Eric Brussel}
\email{ebrussel@calpoly.edu\\mgoertz@calpoly.edu}
\thanks{This research was generously supported by the William and Linda Frost Fund in the Cal Poly Bailey College of Science and Mathematics.}
\author[M. E. Goertz]{Madeleine E. Goertz\\California Polytechnic State University\\San Luis Obispo, USA}

\begin{document}
\small
\title[The Torus of Triangles]
{The Torus of Triangles}

\begin{abstract}
We prove the 2-torus $\T$, an abelian linear algebraic group, is 
a fine moduli space of labeled, oriented, possibly-degenerate inscribable similarity classes of triangles, where a triangle is {\it inscribable} if it can be inscribed in a circle.
A natural action by the dihedral group $D_6$ defines a quotient stack $[\T/D_6]$, which is
the stack of absolute (unlabeled, unoriented) possibly-degenerate inscribable classes.
We show the main triangle types form distinguished algebraic 
substructures: subgroups, cosets,
and elements of small order, and we apply the natural metric on $\T$ to compare them.
\end{abstract}

\subjclass{14C05, 51M05, 60D05}

\maketitle

\tableofcontents

\newpage

\section{\Large Introduction}

In this paper we construct a compactification of the moduli space of labeled, 
oriented similarity classes of triangles in the plane that is 
isomorphic to a 2-torus $\T$, an abelian linear algebraic group. 
A {\it moduli space}, or {\it shape space}, is a continuum of points representing members of a set of 
mathematical objects-of-interest, in this case similarity classes of triangles.
It is useful for studying how objects belong to families, for visualizing 
different types, and for making statistical calculations. 
Often a space of general objects-of-interest is not hard to construct, but
must be compactified using degenerate objects in order to be useful.

Our main result is based on the parameterization of classes
by triples of angles, which naturally restricts the degenerate classes in the compactification 
to those we call {\it inscribable}.
Thus we submit a fine moduli space of 
labeled, oriented, possibly-degenerate {\it inscribable} triangle classes.

One common model of labeled nondegenerate triangle classes, often called the {\it triangle of triangles},
consists of ordered triples of positive (interior) angles adding to $\pi$ in $\R^3$,
see e.g. \cite[Figure 2]{ES15}. 
Our paper was motivated by
the observation that in the triangle-of-triangles model, the triangles in 
certain continuous ``Poncelet'' families inscribed in a circle
appeared to turn over, or reverse orientation, as they crossed the degenerate border.
By incorporating orientation-reversal across this border, we found that 
the triangle of triangles and its oppositely oriented counterpart
glued together to
form a torus $\T$, equipped with a natural
abelian algebraic group structure.
We called $\T$, of course, the {\it torus of triangles}.

Since the torus is a group, we have a pairing, whereby two triangles combine to form a third. 
Consequently our space satisfies the uniformity criterion in \cite{Port},
and addressed in \cite{CNSS}.
Moreover, the standard distinguished triangle types have distinguished group-theoretic significance.
The natural metric on the torus allows us to compute ratios of different types: We find that
the obtuse to acute ratio is $(3:1)$, while the obtuse-isosceles to acute-isosceles ratio is $(1:1)$.
The isosceles to right ratio is $(2\sqrt 3:1)$. 

The shape space $\T$ admits an action by the dihedral group $D_6$ ($12$ elements),
and the orbit space is the (coarse) moduli space of absolute 
(unlabeled, unoriented) possibly-degenerate inscribable classes.
We thus obtain a stack of possibly-degenerate inscribable classes as the quotient $[\T/D_6]$.

\subsubsection{\bf Literature}
Though the torus is a natural and simple compactification of the space
of labeled, oriented, nondegenerate classes, 
we have not found it as such in the literature, where the most commonly named shape space
is a sphere.
The sphere is featured in 
\cite{Kend84}, \cite{Iwai}, \cite{Beh}, \cite{Montgomery}, \cite{ES15}, \cite{CNSS}, and others. 
The principal difference between the torus and the sphere is that the latter is simply connected,
while the former has ``nontrivial loops'', which are continuous families of triangles that can't be evolved into
a single triangle shape. We give an example in Example~\ref{arise}.

The problem of constructing a moduli space of triangles is frequently posed as an elementary
exercise that demonstrates the properties and challenges of more complex moduli spaces,
after a remark attributed to M. Artin. 
Along these lines,
Behrend uses it in \cite{Beh} to introduce the theory of algebraic stacks.
But the idea of constructing a space of triangles goes way back,
and was suggested in the 19C by Lewis Carroll, 
who in ``Pillow Problem 58'' asked
for the probability that a randomly chosen triangle is obtuse (see \cite{Dod58}).
As pointed out in \cite{CNSS}, the first mention of it in the literature appears to be
in {\it The Lady's and Gentleman's Diary} \cite{WSB}, where W.S.B. Woolhouse said:
``In a given circle a regular polygon is inscribed, and lines are drawn from each 
of its angles to the center of the circle. Required the ratio of the number of the 
triangles which are acute-angled to that of those which are obtuse-angled.''

In \cite{Port} Portnoy suggests a moduli space of triangles 
should admit a transitive action by a compact group, in order that
the various classes carry equal priority, and different regions can be assigned
finite measures (since the group is compact). The ``right'' distribution should
align with the spherically symmetrical construction $\P^5(\R)$ using 
the six coordinates of the triangle's three vertices (up to scalar multiplication).
This idea is used in \cite{ES15}.
Since our construction actually {\it is} a compact group, it admits a transitive action
by a compact group, and so passes Portnoy's test, but we include degenerate types
that do not come from $\P^5(\R)$. Some of our measurements
agree with those of \cite{Port}, though our measure itself
appears to be different, as noted in \cite[1.1]{ES15}. 

Our work intersects deep results on compactifications of configuration spaces
of $n$-tuples of points on a manifold or variety, appearing in
\cite{Ful94}, \cite{AxSi94}, and \cite{Sinha}.
These papers treat the problem of representing degenerate configurations of points -- 
which in the $n=3$ case represent degenerate triangles --
essentially as we do: by appending a tangent direction to the multiple point.
In \cite{Ful94} and \cite{AxSi94} this is achieved by blowing up at the degenerate locus; in \cite{Sinha} it is
more explicitly combinatorial and applies more generally to smooth manifolds.
We thank the referee for pointing out these important papers.

\subsubsection*{Acknowledgements}
We thank Sean Gasiorek for sparking their interest
in spaces of triangles at a Cal Poly colloquium on the subject of 
triangle billiards and Poncelet's Porism.

\section{\Large Definitions}

A triangle is the geometric figure obtained by connecting
three points in the plane with straight lines.
In order to make a proper space of triangles we need to accommodate
degenerate examples, which accounts for some of the technicalities in
the following definition.

\begin{Definition}\label{triangles}\rm
A {\it triangle} is a plane figure consisting of three vertices, three straight edges
connecting them, and three  
interior angles. We label the vertices $(A,B,C)\in\C^3$, the side vectors
$(\a,\b,\c)=(C-B,A-C,B-A)\in\C^3$, and
the interior angles $(\alpha,\beta,\gamma)\in(-\pi,\pi]^3$, where
\begin{enumerate}[(a)]
\item\label{(a)}
$(\alpha,\beta,\gamma)\neq(0,0,0)$,
$\alpha+\beta+\gamma\in\{-\pi,0,\pi\}$,
and if $\alpha\beta\gamma\neq 0$ then $\alpha,\beta,\gamma$ have the same sign.
\item\label{(b)}
Let ${\sf Arg}(\bs z)$, the {\it principal argument} of a nonzero complex number $\bs z$,
be the unique real number $\xi$ in the interval $(-\pi,\pi]$ such that $\bs z=|\bs z|e^{\i\xi}$. 
Then \begin{enumerate}[$\bullet$]
\item
$\alpha=\Arg(-\b/\c)$ when $\b\c\neq 0$, 
\item
$\beta=\Arg(-\c/\a)$ when $\c\a\neq 0$, 
\item
$\gamma=\Arg(-\a/\b)$ when $\a\b\neq 0$.
\end{enumerate}
\item
A triangle is {\it degenerate}, or {\it collinear}, if the vertices are collinear, and otherwise {\it nondegenerate}.
\item
A nondegenerate triangle has {\it positive orientation} if the directed graph
defined by its side vectors $(\a,\b,\c)$
goes counterclockwise in the plane, and {\it negative orientation} if it goes clockwise. 
Equivalently, it has positive orientation if its interior angles are positive, and
negative orientation if they are negative.
A degenerate triangle has {\it zero orientation}.
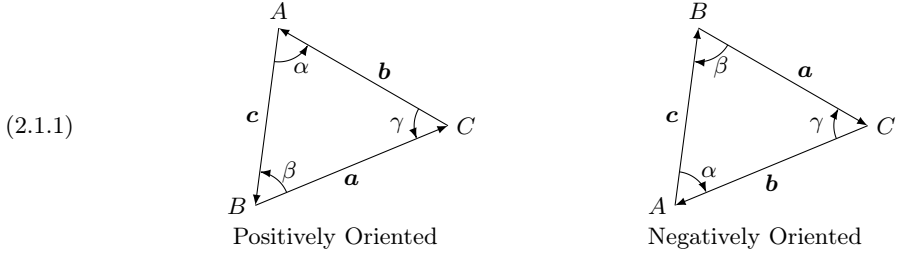
\begin{equation}\label{oriented}
\qquad\qquad\begin{minipage}{.3\textwidth}
\centering
\begin{tikzpicture}[scale=.75]
\coordinate (A) at (-1,{sqrt(3)});
\coordinate (B) at ({-sqrt(2)},{-sqrt(2)});
\coordinate (C) at (2,0);
\draw[->,-latex] (A)--(B);
\draw[->,-latex] (B)--(C);
\draw[->,-latex] (C)--(A);
\node[right] at (C) {$C$};
\node[above] at (A) {$A$};
\node[left] at (B) {$B$};
\node[right] at (1/2+.1,{sqrt(3)/2+.1}) {$\b$};
\node[left] at ({(-sqrt(2)-1)/2},{(sqrt(3)-sqrt(2))/2}) {$\c$};
\node[below] at ({(-sqrt(2)+2)/2},{-sqrt(2)/2}) {$\a$};
\node at ({-sqrt(2)+.6},{-sqrt(2)+.6}) {$\beta$};
\node at ({-1+.4},{sqrt(3)-.7}) {$\alpha$};
\node at (1.1,0) {$\gamma$};
\draw[domain=-1.701696:-.523599,->,-latex] plot ({-1+.6*cos(deg(\x))},{sqrt(3)+.6*sin(deg(\x))}); 
\draw[domain=.392699:1.439897,->,-latex] plot ({-sqrt(2)+.6*cos(deg(\x))},{-sqrt(2)+.6*sin(deg(\x))}); 
\draw[domain=2.617994:3.534292,->,-latex] plot ({2+.6*cos(deg(\x))},{.6*sin(deg(\x))}); 
\node at (0,-2) {Positively Oriented};
\end{tikzpicture}
\end{minipage}
\begin{minipage}{.5\textwidth}
\centering
\begin{tikzpicture}[scale=.75]
\coordinate (B) at (-1,{sqrt(3)});
\coordinate (A) at ({-sqrt(2)},{-sqrt(2)});
\coordinate (C) at (2,0);
\draw[->,-latex] (A)--(B);
\draw[->,-latex] (B)--(C);
\draw[->,-latex] (C)--(A);
\node[right] at (C) {$C$};
\node[left] at (A) {$A$};
\node[above] at (B) {$B$};
\node[right] at (1/2+.1,{sqrt(3)/2+.1}) {$\a$};
\node[left] at ({(-sqrt(2)-1)/2},{(sqrt(3)-sqrt(2))/2}) {$\c$};
\node[below] at ({(-sqrt(2)+2)/2},{-sqrt(2)/2}) {$\b$};
\node at ({-sqrt(2)+.6},{-sqrt(2)+.6}) {$\alpha$};
\node at ({-1+.4},{sqrt(3)-.7}) {$\beta$};
\node at (1.1,0) {$\gamma$};
\draw[domain=-.523599:-1.701696,->,-latex] plot ({-1+.6*cos(deg(\x))},{sqrt(3)+.6*sin(deg(\x))}); 
\draw[domain=1.439897:.392699,->,-latex] plot ({-sqrt(2)+.6*cos(deg(\x))},{-sqrt(2)+.6*sin(deg(\x))}); 
\draw[domain=3.534292:2.617994,->,-latex] plot ({2+.6*cos(deg(\x))},{.6*sin(deg(\x))}); 
\node at (0,-2) {Negatively Oriented};
\end{tikzpicture}
\end{minipage}
\end{equation}
\item\label{(d)}
A triangle is {\it inscribable} if it can be inscribed in a circle,
and when all vertices are identified then $\alpha\beta\gamma=0$.
\item\label{(e)}
Two labeled, oriented, possibly-degenerate inscribable triangles are {\it similar} if their triples
of vertices differ by a nonzero complex number and their interior angles are equal $\pmod\pi$.
An (inscribable) {\it triangle class}, or just {\it class}, is the resulting equivalence class.
\end{enumerate}
\noindent Write 
\[\Delta=\{\left((A,B,C);(\alpha,\beta,\gamma)\right)\}\subset\C^3\times(-\pi,\pi]^3\] 
for the set of all labeled, oriented, possibly-degenerate triangles,
$\Delta_{\rm insc}$ for the inscribable subset,
and $[\Delta_{\rm insc}]$ for the set of inscribable classes.
\end{Definition}

\begin{Remarks}
The definition of similarity in \eqref{(e)} for inscribable triangles 
agrees with the standard one when all vertices are distinct,
but admits a range of (degenerate) classes where two or three vertices coincide.
These extra classes, which we call {\it double} and {\it triple points}, 
are forced by our compactification construction.
By \eqref{(a)} and \eqref{(b)}, a double point class with $A=C\neq B$ has the form
$[(A,B,A);(\alpha,0,-\alpha)]$, where $\alpha$ is any value in $(-\pi,\pi]$.
We will show they are meaningfully depicted with a direction at the doubled vertex:
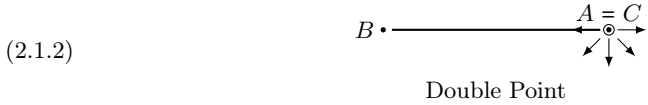
\begin{equation}\label{doublepoints}
\begin{minipage}{.35\textwidth}
\begin{tikzpicture}
\node at (0,0) (A) {};
\node at (3,0) (B) {};
\draw[thick] (A) -- (B);
\draw[->,-latex] (B)--({3+.5*cos(180)},{.5*sin(180)});
\draw[->,-latex] (B)--({3+.5*cos(225)},{.5*sin(225)});
\draw[->,-latex] (B)--({3+.5*cos(270)},{.5*sin(270)});
\draw[->,-latex] (B)--({3+.5*cos(315)},{.5*sin(315)});
\draw[->,-latex] (B)--({3+.5*cos(0)},{.5*sin(0)});
\fill (0,0) circle (1pt);
\fill (3,0) circle (1pt);
\draw (3,0) circle (2pt);
\node[left] at (A) {$B$};
\node[above] at (B) {$A=C$};
\node at (1.5,-.8) {Double Point};
\end{tikzpicture}
\end{minipage}
\end{equation}
\noindent
Triple point classes like
$[(A,A,A);(\alpha,\beta,\gamma)]$ with $\alpha\beta\gamma=0$ in \eqref{(d)}
are limits of sequences in $[\Delta_{\rm insc}]$.
For example, $((A,A,A);(\pi,0,0))$ is the limit as $\ep$ goes to $0$ 
of the sequence with vertices $(e^{\ep\i},e^{2\ep\i},1)$.
Similarly, the limit as $B$ approaches $A=C$ above gives a point
$((A,A,A);(\alpha,0,-\alpha))$ in $\Delta_{\rm insc}$,
a ``tripled double point''.
\end{Remarks}

\section{\Large The Torus}\label{thetorus}
First we construct the space that will parameterize the inscribable classes.
The quotient map 
\[\R^3\lr(\R/\pi)^3\equiv\R/\pi\Z\times\R/\pi\Z\times\R/\pi\Z\] 
is a smooth covering of abelian Lie groups by
\cite[Example 21.14]{Lee}.
In particular $(\R/\pi)^3$ is a smooth manifold, the $3$-torus.
Let
$\T\subset(\R/\pi)^3$ be the image 
of the linear subspace $L:X+Y+Z=0$ of $\R^3$. Then $\T$ is a Lie subgroup by \cite[Theorem 21.27]{Lee},
with the induced pointwise additive composition law.
Taking representatives of the coordinates on $L$ in the fundamental domain $[0,\pi)$
allows us to visualize $\T$ as the intersection of a solid cube of side-length $\pi$
with the parallel planes $\frac k3(\pi,\pi,\pi)+L$ ($k=0,1,2$):
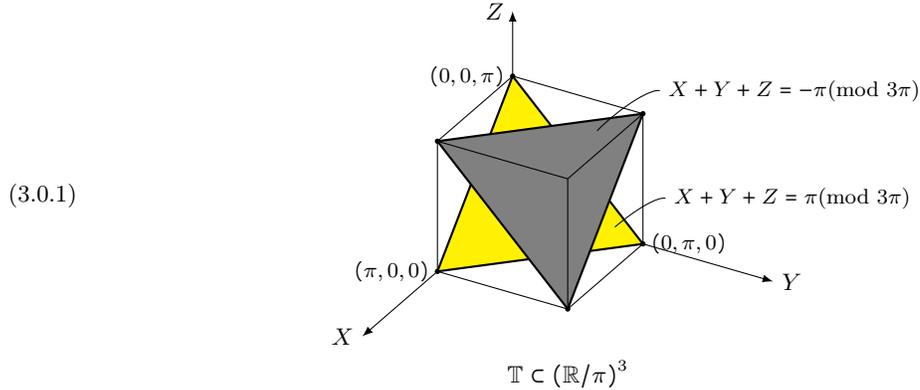
\begin{equation}\label{cube}
\begin{minipage}{.4\textwidth}
\tdplotsetmaincoords{60}{120}
\begin{tikzpicture}[scale=2,tdplot_main_coords]
\coordinate (A) at (0,0,0);
\coordinate (B) at (1,0,0);
\coordinate (C) at (1,1,0);
\coordinate (D) at (0,1,0);
\coordinate (E) at (0,0,1);
\coordinate (F) at (0,1,1);
\coordinate (G) at (1,0,1);
\coordinate (H) at (1,1,1);
\draw (A)--(B)--(C)--(D)--cycle; 
\draw (A)--(E)--(F)--(D)--cycle; 
\draw (A)--(B)--(G)--(E)--cycle; 
\draw[->,-latex] (A)--(2,0,0) node[left] {$X$};
\draw[->,-latex] (A)--(0,2,0) node[right] {$Y$};
\draw[->,-latex] (A)--(0,0,1.5) node[left] {$Z$};
\fill (A) circle (.5pt);
\fill (B) circle (.5pt) node[left] {\footnotesize $(\pi,0,0)$};
\fill (D) circle (.5pt) node[right] {\footnotesize $(0,\pi,0)$};
\fill (C) circle (.5pt);
\fill (E) circle (.5pt) node[left] {\footnotesize $(0,0,\pi)$};
\fill (F) circle (.5pt);
\fill (G) circle (.5pt);
\fill (H) circle (.5pt);
\fill[yellow, opacity=0.5] (1,0,0) -- (0,1,0) -- (0,0,1) -- cycle;
\draw[thick] (B)--(D)--(E)-- cycle;
\fill[gray, opacity=0.5] (1,1,0) -- (1,0,1) -- (0,1,1) -- cycle;
\draw[thick] (C)--(G)--(F)-- cycle;
\draw (G)--(H)--(F); 
\draw (H)--(C); 
\node at (1,1,-.5) {$\T\subset(\R/\pi)^3$};
\draw (.35,.85,1) arc (290:180:.5 and .1) node[right] {\footnotesize $X+Y+Z=-\pi\pmod{3\pi}$};
\draw (.2,.9,.2) arc (270:180:.5 and .1) node[right] {\footnotesize $X+Y+Z=\pi\pmod{3\pi}$};
\end{tikzpicture}
\end{minipage}
\end{equation}
The boundaries of the triangular regions, where one angle is
$0\pmod\pi$, are identified
on opposite faces of the cube, and all 8 vertices, where all
angles are $0\pmod\pi$, are identified to a single point.
As a set we have
\[\T=\left\{(\alpha,\beta,\gamma)\in(\R/\pi)^3:\alpha+\beta+\gamma=0\pmod\pi\right\}\]

\begin{Lemma}\label{torus}
$\T\subset(\R/\pi)^3$ is naturally isomorphic to a $2$-torus.
\end{Lemma}

\begin{proof}
The identifications of $(\R/\pi)^3$ shows the two triangles of $\T$ form a standard 
torus polygon:
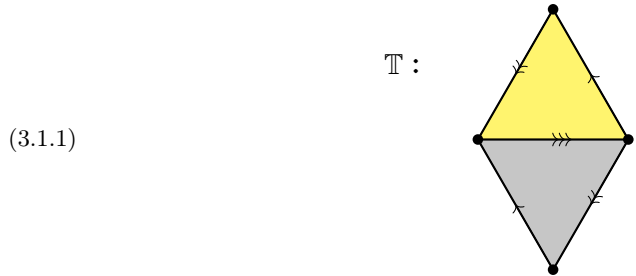
\begin{equation}\label{polygon}
\begin{minipage}{.3\textwidth}
\begin{tikzpicture}[scale=2]
\coordinate (A) at (-1/2,0);
\coordinate (B) at (0,{sqrt(3)/2});
\coordinate (C) at (1/2,0);
\coordinate (D) at (0,{-sqrt(3)/2});
\draw[thick] (A)--(B)--(C)--(D)--(A)--(C);
\fill (A) circle (1pt);
\fill (B) circle (1pt);
\fill (C) circle (1pt);
\fill (D) circle (1pt);
\draw[->>] (B)--(-1/4,{sqrt(3)/4});
\draw[->>] (C)--(1/4,{-sqrt(3)/4});
\draw[->>>] (A)--(1/8,0);
\draw[->] (C)--(1/4,{sqrt(3)/4});
\draw[->] (D)--(-1/4,{-sqrt(3)/4});
\begin{scope}[on background layer]
\draw[fill=darkgray!30] (A)--(D)--(C)--(A);
\draw[fill=yellow!70] (A)--(B)--(C)--(A);
\end{scope}
\node at (-1,.5) {\large $\T:$};
\end{tikzpicture}
\end{minipage}
\end{equation}
\end{proof}
\begin{Corollary}\label{algvar}
$\T$ is an abelian linear algebraic group.
In particular it is an algebraic variety, and contains closed subgroups 
corresponding to the linear subspaces in \eqref{cube}.
\end{Corollary}

\begin{proof}
Since $\T$ is a compact Lie group, it is a linear algebraic group by \cite[Ch.5, 2.5, Theorem 12]{OV}.
By the corollary to that theorem there is a one-to-one correspondence
between Lie subgroups and closed algebraic subgroups, and we
obtain the second statement since every linear subspace of 
$\R^3$ is a Lie subgroup.
\end{proof}

\section{\Large Parameterizing Inscribable Classes}

The following theorem effectively
shows $\T$ is the moduli space, or shape space, of 
labeled, oriented, possibly degenerate inscribable triangle classes.

\begin{Theorem}\label{map}
$\T$ is naturally isomorphic to the set $[\Delta_{\rm insc}]$ of labeled, oriented, possibly-degenerate
inscribable classes.
More precisely, the natural map 
\begin{align*}
\rho:[\Delta_{\rm insc}]&\lr\T\subset(\R/\pi)^3\\
[\left((A,B,C);(\alpha,\beta,\gamma)\right)]&\lm(\alpha,\beta,\gamma)\pmod\pi
\end{align*} 
is a bijection, taking positively oriented classes to the subset $X+Y+Z=\pi\pmod{3\pi}$ of $(\R/\pi)^3$,
negatively oriented classes to $X+Y+Z=-\pi\pmod{3\pi}$,
and degenerate classes to the boundary $\alpha\beta\gamma=0\pmod\pi$ in \eqref{cube}.
\end{Theorem}

\begin{proof}
The theorem says that a labeled, oriented, possibly-degenerate 
inscribable class is completely determined
by its triple of angles $\pmod\pi$; that every triple
$\pmod\pi$ arises from some class; and that the two shaded triangular regions in
\eqref{cube} correspond to the two possible orientations.

All three angles of a nondegenerate class belong to either $(0,\pi)$ or $(-\pi,0)$,
depending on orientation, so its image is in the interior $\alpha\beta\gamma\neq 0$ of $\T$,
and is uniquely determined by its angle values $\pmod\pi$ together with its orientation.
If the class is positively oriented then $\alpha+\beta+\gamma=\pi\pmod{3\pi}$,
so the image is in the yellow region in \eqref{cube},
and if negatively oriented $\alpha+\beta+\gamma=-\pi\pmod{3\pi}$, in the
gray region. 
Therefore $\rho$ distinguishes orientation, and is injective on nondegenerate classes.

A degenerate class with $\gamma=0$ has form $[(A,A,C);(\alpha,-\alpha,0)]$ with
$\alpha\in(-\pi,\pi]$, and maps to the identified faces of the cube \eqref{cube} parallel
to the $XY$-plane. Suppose
$\rho([(A,A,C);(\alpha,-\alpha,0)])=\rho([(A',A',C');(\alpha',-\alpha',0)])$. 
Since there are only two, the vertices automatically differ by a nonzero complex constant. Since
$\alpha'=\alpha\pm\pi$, the interior angles are equal $\pmod\pi$.
Therefore the triangles are similar
by Definition~\ref{triangles}\eqref{(e)}.
Clearly
this applies to the double point degenerates $\alpha=0$ and $\beta=0$ as well,
which go to the faces of the cube parallel to the $YZ$ and $XZ$-planes, respectively.
Thus $\rho$ remains injective if we include these classes.
The solitary triple point degenerate class is $[(1,1,1);(\pi,0,0)]$ maps to $(0,0,0)\pmod\pi$,
and since this point has not been taken, we conclude $\rho$ is injective in general.

If $\alpha+\beta+\gamma=\pi\pmod{3\pi}$ and $\alpha\beta\gamma\neq 0\pmod\pi$ in $\T$
then there exist real representatives in $(0,\pi)$ adding to $\pi$, and these are
the angles of a positively oriented triangle, by Euclid. Similarly if the sum is $-\pi$
we find a negatively oriented triangle. Therefore 
$\rho$ is surjective on the interior regions of $\T$, dividing positively and negatively oriented
classes as above.
If $\gamma=0$ and $\alpha\neq 0\pmod\pi$ then $(\alpha,\beta,\gamma)\in\T$ is in the image of
the class $[(A,A,C);(\alpha,-\alpha,0)]$, and similarly if $\alpha=0$ or $\beta=0$.
We conclude $\rho$ is surjective, hence bijective, and divides classes as claimed.
\end{proof}

\begin{Remark}
We prove more rigorously in Corollary~\ref{moduli} 
that $\T$ is a fine moduli space; of labeled, oriented, possibly-degenerate
inscribable classes. As the closure of the interior $\alpha\beta\gamma\neq 0$, it then represents 
a compactification of the space of labeled, oriented 
nondegenerate classes.
\end{Remark}

\subsection{Families and Degenerates}
We next invert the correspondence $\rho$ in Theorem~\ref{map} to give an explicit family
of inscribed triangles corresponding to the points of $\T$. This will help explain
how the degenerate classes in $[\Delta_{\rm insc}]$, which are
classes with one or more zero interior angle, arise naturally as limit points of sequences of
nondegenerates.
First, a lemma from Euclid.

\begin{Lemma}\label{chord}
For points $B$ and $C$ on a circle, let $\Arc(BC)$ be the {\it oriented} (counterclockwise)
arc from $B$ to $C$, in radians.
Suppose $((A,B,C);(\alpha,\beta,\gamma))$ is nondegenerate and inscribed in a circle.
Then $\alpha=\frac 12\Arc(BC)\pmod\pi$, 
and $\alpha$ equals the oriented angle from the tangent line at $B$ to 
the chord $BC$, or equivalently from $BC$ to the tangent at $C$.
\end{Lemma}

\begin{proof}
The positively oriented case $\alpha=\frac 12\Arc(BC)$ is
Proposition 20 of Euclid's Book III, which states,
{\it In a circle the angle at the center is double the angle at 
the circumference when the angles have the same circumference as base} (\cite{Euclid}). 
For negatively oriented $\alpha'<0$ we have $\alpha'=-\frac 12\Arc(CB)=\frac 12(\Arc(BC)-2\pi)=\frac 12(\Arc(BC)-\pi)$,
so again $\alpha'=\frac 12\Arc(BC)\pmod\pi$, as desired. 

Diagram \eqref{alpha} shows how to compute the opposing angle from the tangent at $B$
to the chord $BC$ using similar triangles $MOB\sim MBI$. The last statement is immediate.
\begin{equation}
\begin{minipage}{.4\textwidth}
\begin{tikzpicture}[scale=1.25]\label{alpha}
\coordinate (O) at (0,0) node[left] at (O) {$O$};
\coordinate (B) at (0,-1) node[below] at (B) {$B$};
\coordinate (C) at ({1/sqrt(2)},{1/sqrt(2)}) node[right] at (C) {$C$};
\coordinate (M) at ({1/(2*sqrt(2))},{(-sqrt(2)+1)/(2*sqrt(2))}) node[right] at (M) {$M$};
\coordinate (I) at ({(2+sqrt(2))/sqrt(2)},-1) node[right] at (I) {$I$};
\coordinate (A) at (-1/2,{sqrt(3)/2}) node[left] at (A) {$A$};
\coordinate (A') at (1,0) node[right] at (A') {$A'$};
\draw (O) circle (1);
\draw (O)--(B)--(C)--(O)--(I)--(C);
\draw (I)--(-.6,-1);
\draw[dashed] (B)--(A)--(C);
\draw[dashed] (B)--(A')--(C);
\fill (O) circle (.75pt);
\fill (B) circle (.75pt);
\fill (C) circle (.75pt);
\fill (M) circle (.75pt);
\fill (I) circle (.75pt);
\fill (A) circle (.75pt);
\fill (A') circle (.75pt);
\draw[domain=0:1.178097,->,-latex] plot ({.25*cos(deg(\x))},{-1+.25*sin(deg(\x))}); 
\draw[domain=-3.141593:-4.31969,->,-latex] plot ({.15*cos(deg(\x))},{-1+.15*sin(deg(\x))}); 
\draw[domain=-1.570796:-.392699,->,-latex] plot ({.25*cos(deg(\x))},{.25*sin(deg(\x))});
\draw[domain=-1.308997:-.1309,->,-latex] plot ({-1/2+.25*cos(deg(\x))},{(sqrt(3)/2)+.25*sin(deg(\x))});
\draw[domain=-135:-247.5,->,-latex] plot ({1+.15*cos(\x)},{.15*sin(\x)}); 
\draw (.85,0) arc (-135:-260:.8 and .8/3) node[right] {\footnotesize $\alpha'=\alpha-\pi$}; 
\node at (.2,-.3) {\footnotesize $\alpha$};
\node at (.35,-.85) {\footnotesize $\alpha$};
\node at (-.2,.65) {\footnotesize $\alpha$};
\node[left] at (-.1,-.85) {\footnotesize $\alpha-\pi$};
\end{tikzpicture}
\end{minipage}
\end{equation}
\end{proof}

The next result gives an inscribed family over $\T$.

\begin{Theorem}\label{family}
The map $\rho$ of Theorem~\ref{map} has inverse
\begin{align*}
\sigma:\T&\lr[\Delta_{\rm insc}]\\
(\alpha,\beta,\gamma)\pmod\pi&\lm
\left[(e^{2\beta\i},e^{-2\alpha\i},1);(\alpha,\beta,\gamma)\right]
\quad\text{if $(\alpha,\beta,\gamma)\neq (0,0,0)\pmod\pi$}\\
(0,0,0)\pmod\pi&\lm\left[(1,1,1);(\pi,0,0)\right]
\end{align*}
where the interior angles in the image are real representatives in $(-\pi,\pi]$ 
adding to $\pm\pi$, of the same sign when $\alpha\beta\gamma\neq 0$.
\end{Theorem}

\begin{proof}
The map $\sigma$ as stated clearly inverts $\rho$ if it is well-defined, i.e.,
if real representatives $(\alpha,\beta,\gamma)$ can be chosen as stated, and
then if the resulting vertex-angle pairings are actually triangle classes. 
Since the choice of interior angles exists and is unique by Figure~\ref{cube}, 
it remains to check the configuration
of vertices is compatible with the interior angles.

Assume $\alpha\beta\gamma\neq 0$.
Then $\sigma$ is illustrated below in Figure~\ref{euclid}, with the oriented arcs computed using 
the vertices' arguments:
\begin{equation}\label{euclid}
\begin{minipage}{.7\textwidth}
\begin{tikzpicture}[scale=1.15]
\coordinate (A) at ({cos(110)},{sin(110)});
\coordinate (B) at ({cos(210)},{sin(210)});
\coordinate (C) at (1,0);
\draw (0,0) circle (1);
\filldraw (1,0) circle (.5pt);
\filldraw (A) circle (.5pt);
\filldraw (B) circle (.5pt);
\filldraw (0,0) circle (.5pt);
\draw (A)--(B)--(C)--(A);
\draw (1,0)--(0,0)--(A);
\draw (0,0)--(B);
\draw[domain=10:110,->,-latex] plot ({1.15*cos(\x)},{1.15*sin(\x)}); 
\draw[domain=253:328,->,-latex] plot ({cos(110)+.25*cos(\x)},{sin(110)+.25*sin(\x)}); 
\draw[domain=15:70,->,-latex] plot ({cos(210)+.25*cos(\x)},{sin(210)+.25*sin(\x)}); 
\draw[domain=-150:-10,->,-latex] plot ({1.2*cos(\x)},{1.2*sin(\x)}); 
\node[right] at (.9,.9) {\footnotesize $2\beta$};
\node[right] at (.5,-1.2) {\footnotesize $2\alpha$};
\node at ({cos(210)+.35},{sin(210)+.3}) {$\beta$};
\node at ({cos(110)+.15},{sin(110)-.35}) {$\alpha$};
\node[right] at (C) {$C=1$};
\node[left] at (B) {$B=e^{-2\alpha\i}$};
\node[above left] at (A) {$A=e^{2\beta\i}$}; 
\coordinate (D) at ({4.2+cos(110)},{sin(110)});
\coordinate (E) at ({4.2+cos(210)},{sin(210)});
\coordinate (F) at (5.2,0);
\draw (4.2,0) circle (1);
\filldraw (5.2,0) circle (.5pt);
\filldraw (D) circle (.5pt);
\filldraw (E) circle (.5pt);
\filldraw (4.2,0) circle (.5pt);
\draw (D)--(E)--(F)--(D);
\draw (5.2,0)--(4.2,0)--(D);
\draw (4.2,0)--(E);
\draw[domain=10:110,->,-latex] plot ({4.2+1.15*cos(\x)},{1.15*sin(\x)}); 
\draw[domain=253:328,->,-latex] plot ({4.2+cos(110)+.25*cos(\x)},{sin(110)+.25*sin(\x)}); 
\draw[domain=15:70,->,-latex] plot ({4.2+cos(210)+.25*cos(\x)},{sin(210)+.25*sin(\x)}); 
\draw[domain=-150:-10,->,-latex] plot ({4.2+1.2*cos(\x)},{1.2*sin(\x)}); 
\node[right] at (5.1,.9) {\footnotesize $-2\alpha$};
\node[right] at (4.7,-1.2) {\footnotesize $-2\beta$};
\node at ({4.2+cos(210)+.35},{sin(210)+.3}) {$-\alpha$};
\node at ({4.2+cos(110)+.15},{sin(110)-.35}) {$-\beta$};
\node[right] at (F) {$C=1$};
\node[left] at (E) {$A=e^{2\beta\i}$};
\node[above left] at (D) {$B=e^{-2\alpha\i}$};
\node at (2,-1.8) {Relation Between Interior Angles and Arcs};
\end{tikzpicture}
\end{minipage}
\end{equation}
The interior angles and corresponding arcs in \eqref{euclid} are compatible by 
Lemma~\ref{chord}, hence the angles appear in exponents at the given vertices,
and the map is well-defined on $\alpha\beta\gamma\neq 0$.

Now suppose $\alpha\beta\gamma=0$. 
By Definition~\ref{triangles}\eqref{(a)} we know at least one angle is nonzero. 
If $(\alpha,\beta,0)\in\T$ and $\alpha\neq 0\pmod\pi$, then $\beta=-\alpha$, so
\[\sigma(\alpha,\beta,\gamma)=[(e^{-2\alpha\i},e^{-2\alpha\i},1);(\alpha,-\alpha,0)]\] which is well-defined.
Similarly if $\alpha=0$ and $\beta\neq 0\pmod\pi$, or $\beta=0$ and $\alpha\neq 0\pmod\pi$, we find the image
of $\sigma$ is $[(e^{2\beta\i},1,1);(0,\beta,-\beta)]$ and $[(1,e^{-2\alpha\i},1);(\alpha,0,-\alpha)]$,
respectively, both well-defined. Finally the remaining point $(0,0,0)\in\T$ maps to $[(1,1,1);(\pi,0,0)]$,
as desired. 
We conclude $\sigma$ is well-defined, hence it inverts $\rho$.
\end{proof}

\begin{Remark}
The map $\sigma$ gives inscribed triangles over $\T$
naturally representing every class in $[\Delta_{\rm insc}]$.
The existence of such a family is an important feature of a (fine) moduli space,
a sort of testimonial to the fine-ness.
Though $\T$ itself tautologically provides a family of classes via $\rho$,
Theorem~\ref{family} goes further by specifying an explicit set of triangle representatives. 
We define families more rigorously in Section~\ref{moduliquotient}.
\end{Remark}

\subsection{A Noncontractible Loop of Classes}
\begin{Example}\label{arise}
We use $\sigma$ to explicitly write down a noncontractible loop of triangle classes,
and at the same time show how the degenerate double point classes arise in inscribed families,
in this case defined by a fixed chord $BC$.

Fix $\alpha_0\in(0,\pi)$. Then the set 
$\{(\alpha_0,\beta,-(\alpha_0+\beta))\pmod\pi\in\T\}$
is a ``circle'' of points on $\T$, both positively and negatively
oriented, parameterized by $\beta\in\R/\pi$:
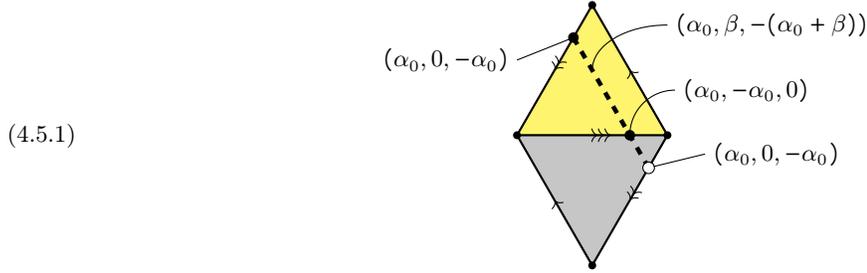
\begin{equation}\label{fixedalpha}
\begin{minipage}{.3\textwidth}
\begin{tikzpicture}[scale=2]
\coordinate (A) at (-1/2,0); 
\coordinate (B) at (0,{sqrt(3)/2}); 
\coordinate (C) at (1/2,0); 
\coordinate (D) at (0,{-sqrt(3)/2}); 
\coordinate (E) at (-1/8,{3*sqrt(3)/8}); 
\coordinate (F) at (3/8,{-sqrt(3)/8}); 
\coordinate (G) at (1/4,0);
\draw[thick] (A)--(B)--(C)--(D)--(A)--(C);
\draw[line width=0.6mm, dashed] (E)--(F);
\fill (A) circle (.75pt);
\fill (B) circle (.75pt);
\fill (C) circle (.75pt);
\fill (D) circle (.75pt);
\fill (E) circle (1pt);
\fill (G) circle (1pt);
\draw (0,{sqrt(3)/4}) arc (180:90:.5 and .3) node[right] {$(\alpha_0,\beta,-(\alpha_0+\beta))$};
\draw (G) arc (180:90:.3 and .3) node[right] {$(\alpha_0,-\alpha_0,0)$};
\draw (E) -- (-1/2,1/2) node[left] {$(\alpha_0,0,-\alpha_0)$};
\draw (F) -- (3/4,-1/8) node[right] {$(\alpha_0,0,-\alpha_0)$};
\draw[->>] (B)--(-1/4,{sqrt(3)/4});
\draw[->>] (C)--(1/4,{-sqrt(3)/4});
\draw[->>>] (A)--(1/8,0);
\draw[->] (C)--(1/4,{sqrt(3)/4});
\draw[->] (D)--(-1/4,{-sqrt(3)/4});
\begin{scope}[on background layer]
\draw[fill=darkgray!30] (A)--(D)--(C)--(A);
\draw[fill=yellow!70] (A)--(B)--(C)--(A);
\end{scope}
\fill[white] (F) circle (1pt); 
\draw (F) circle (1.1pt);
\end{tikzpicture}
\end{minipage}
\end{equation}
This circle of triangles is a noncontractible loop on the torus, which means that the 
triangles cannot be evolved continuously to a single triangle shape.
On a torus embedded in $\R^3$ it could be a loop around the tubular part.
To see what triangles we are dealing with we apply $\sigma$ to obtain
the family
\begin{align*}
\{[(e^{2\beta\i},e^{-2\alpha_0\i},1);(\alpha_0-\pi,\beta,-(\alpha_0+\beta))]
:\beta\in[-\alpha_0,0]\}
\subset[\Delta_{\rm insc}]&\quad\text{if negatively oriented}\\
\{[(e^{2\beta\i},e^{-2\alpha_0\i},1);(\alpha_0,\beta,\pi-(\alpha_0+\beta))]
:\beta\in[0,\pi-\alpha_0]\}
\subset[\Delta_{\rm insc}]&\quad\text{if positively oriented}
\end{align*}
The chord $BC$ is fixed,
and various $\beta$ locate $A$ anywhere on the unit circle:
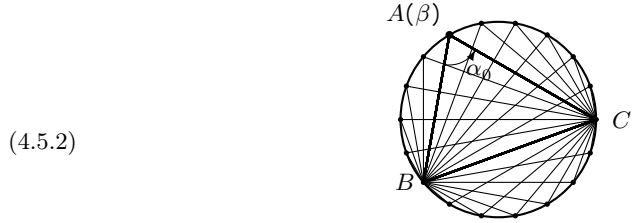
\begin{equation}\label{family1}
\begin{minipage}{.3\textwidth}
\begin{tikzpicture}[scale=1]
\def\radius{1.3}
\draw[thick] (0,0) circle (\radius cm);
\coordinate (target) at ({\radius*cos(-140)},{\radius*sin(-140)}); 
\coordinate (target2) at ({\radius*cos(0)},{\radius*sin(0)}); 
\foreach \angle in {-120, -100, ..., 200}{
\coordinate (vertex) at ({\radius*cos(\angle)},{\radius*sin(\angle)});
\fill (vertex) circle (1pt);
\draw (vertex) -- (target);
\draw (vertex) -- (target2);
\draw[thick] (target)--(target2)--({\radius*cos(120)},{\radius*sin(120)})--(target);
}
\fill (target) circle (1pt);
\fill ({\radius*cos(120)},{\radius*sin(120)}) circle (1.5pt);
\draw[domain=-1.745329:-.523599,->,-latex] plot ({-\radius/2+.4*cos(deg(\x))},{\radius*sqrt(3)/2+.4*sin(deg(\x))});
\node at ({-\radius/2+.4},{\radius*sqrt(3)/2-.5}) {$\alpha_0$};
\node[ right] at ({.1+\radius},0) {$C$};
\node[ left] at (target) { $B$};
\node[above left] at ({\radius*cos(120)},{\radius*sin(120)}) { $A(\beta)$};
\node at (0,-2) {Fixed $\alpha=\alpha_0\pmod\pi$};
\end{tikzpicture}
\end{minipage}
\end{equation}
The values
$\beta=0$ and $\beta=-\alpha_0\pmod\pi$ field the degenerate double points
\begin{align*}
[(1,e^{-2\alpha_0\i},1);(\alpha_0,0,-\alpha_0)]&=[(C,B,C);(\alpha_0,0,-\alpha_0)]\\
[(e^{-2\alpha_0\i},e^{-2\alpha_0\i},1);(\alpha_0,\pi-\alpha_0,0)]&=[(B,B,C);(\alpha_0,\pi-\alpha_0,0)]
\end{align*}
marked on Figure~\ref{fixedalpha}.
Both are represented by the segment $BC$,
but with a different doubled vertex, either $A=B$ or $A=C$.
Since $\alpha_0$ can be any value in $(0,\pi)$, in this way each degenerate inscribable class 
corresponds to a specific chord on the unit circle
with a double point at one end or the other, appearing in a continuous
family of nondegenerate inscribable classes.
By the second part of Lemma~\ref{chord} the degenerate inscribable class
$[(C,B,C);(\alpha_0,0,-\alpha_0)]$ above can be drawn as
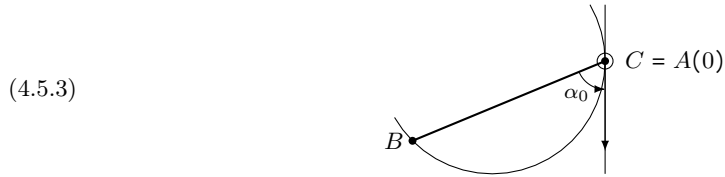
\begin{equation}\label{deginsc}
\begin{minipage}{.3\textwidth}
\begin{tikzpicture}[scale=1.5]
\coordinate (B) at ({-sqrt(2)/2},{-sqrt(2)/2}) node[left] at (B) {$B$};
\coordinate (C) at (1,0);
\node[right] at (1.1,0) {$C=A(0)$};
\fill (B) circle (1pt);
\fill (C) circle (1pt);
\draw (C) circle (2pt);
\draw[thick] (B)--(C);
\draw[->,-latex] (C)--(1,-.8);
\draw (1,-1)--(1,.5);
\draw[domain=-150:30] plot ({cos(\x)},{sin(\x)});
\draw[domain=-2.748894:-1.570796,->,-latex] plot ({1+.25*cos(deg(\x))},{.25*sin(deg(\x))});
\node at (.75,-.3) {\footnotesize $\alpha_0$};
\end{tikzpicture}
\end{minipage}
\end{equation}
Different degenerate inscribable classes
can be distinguished geometrically using directions from the double point:
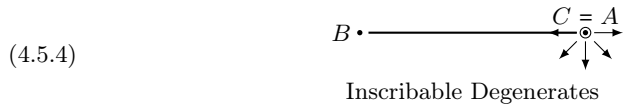
\begin{equation}\label{pixie}
\begin{minipage}{.4\textwidth}
\begin{tikzpicture}
\node at (0,0) (A) {};
\node at (3,0) (B) {};
\draw[thick] (A) -- (B);
\draw[->,-latex] (B)--({3+.5*cos(180)},{.5*sin(180)});
\draw[->,-latex] (B)--({3+.5*cos(225)},{.5*sin(225)});
\draw[->,-latex] (B)--({3+.5*cos(270)},{.5*sin(270)});
\draw[->,-latex] (B)--({3+.5*cos(315)},{.5*sin(315)});
\draw[->,-latex] (B)--({3+.5*cos(0)},{.5*sin(0)});
\fill (0,0) circle (1pt);
\fill (3,0) circle (1pt);
\draw (3,0) circle (2pt);
\node[left] at (A) {$B$};
\node[above] at (B) {$C=A$};
\node at (1.5,-.8) {Inscribable Degenerates};
\end{tikzpicture}
\end{minipage}
\end{equation}
\end{Example}

The point $(0,0,0)\in\T$
uniquely corresponds to the triple point degenerate class $[(1,1,1);(\pi,0,0)]$.
Thus we distinguish three degenerate class types:
\begin{equation}
\begin{minipage}{.6\textwidth}
\begin{tikzpicture}\label{degenerates} 
\filldraw (-2,0) circle (1pt);
\draw (-2,0) circle (2pt);
\draw (-2,0) circle (3pt);
\filldraw (1,0) circle (1pt);
\draw (1,0) circle (2pt);
\filldraw (-1,0) circle (1pt);
\filldraw (4,0) circle (1pt);
\draw (4,0) circle (2pt);
\filldraw (2.5,{sqrt(3)/2}) circle (1pt);
\draw (-3,0) circle (1cm);
\draw (0,0) circle (1cm);
\draw (3,0) circle (1cm);
\draw (-1,0)--(1,0);
\draw (4,0)--(2.5,{sqrt(3)/2});
\node at (-3,-1.5) {triple point};
\node at (0,-1.5) {right double point};
\node at (3.5,-1.5) {general double point};
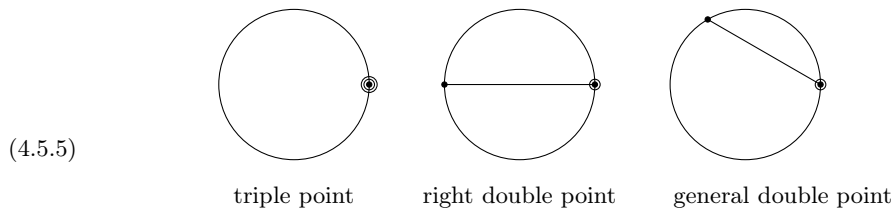
\node at (0,-2.5) {Degenerate Inscribable Triangle Types};
\end{tikzpicture}
\end{minipage}
\end{equation}

\subsection{Missing Degenerates and the Sphere}\label{missing}
The point $(0,0,0)\pmod\pi\in\T$ conceals
the entire set of (degenerate) collinear {\it un}inscribable classes, whose angles are all $0\pmod\pi$, and 
which are classified by the ratios of their (nonzero) 
side lengths $|\a|,|\b|,|\c|$:
\begin{equation}\label{simplept}
\begin{minipage}{.3\textwidth}
\begin{tikzpicture}
\node at (0,0) (B) {};
\node at (3,0) (C) {};
\node at (2,0) (A) {};
\draw[thick] (B)--(A) -- (C);
\fill (A) circle (1pt);
\fill (B) circle (1pt);
\fill (C) circle (1pt);
\node[below] at (A) {$A$};
\node[below] at (B) {$B$};
\node[below] at (C) {$C$};
\node at (1.5,-.8) {Collinear Degenerates};
\end{tikzpicture}
\end{minipage}
\end{equation}
Whenever the shape space of triangle classes is said to be a sphere,
as in \cite{Kend84}, \cite{Montgomery}, \cite{Beh}, and \cite{ES15},
the degenerate classes are of this type.
They are used in \cite[1.1.10]{Beh} to construct
a different compactification
of labeled, oriented, nondegenerate classes, called a {\it bipyramid}. The construction is similar to ours,
describing two triangular regions and gluing them along the degenerate (collinear) boundaries.
The bipyramid is a topological sphere, 
(\cite[Figure 1.21]{Beh}), later recast as the Riemann sphere.

Since neither the torus nor the sphere incorporates all 
degenerate configurations of three vertices, neither gives a complete
 compactification of 
the space of labeled, oriented, nondegenerate triangle classes.
And there are more:
We have also omitted triple points of the form $((1,1,1);(\alpha,\beta,\gamma))$
with $\alpha\beta\gamma\neq 0$, which are realized as limits of the 
nondegenerate triangles in a fixed class as the side vectors are scaled to zero while the angles are held
constant.
We exclude these because they do not arise as limit points in sequences of
nondegenerate inscribable triangles.

In \cite{BGGL} we show that a compactification exists, called {\it Dyck's surface}, 
which takes all ``first order'' degenerate configurations into account,
and projects to both the sphere and the torus.
More precisely, it is the blowup of the sphere at three points 
(``pinched points'' in \cite{Beh}), or the torus at a single point (the degenerate triple point),
and thus it subsumes both constructions.

\section{\Large Group-Theoretic Interpretation of Triangle Types}\label{disting}

In this section we show
the group structure of $\T$
is compatible with triangle types, in the sense that the latter form
basic algebraic structures: elements of finite order, subgroups, and cosets.
By construction (Section~\ref{thetorus})
the composition law is addition of triples of angles $\pmod\pi$.
The identity is $(0,0,0)$, the degenerate triple point.
The main triangle types are illustrated in Figure~\ref{diamond}.
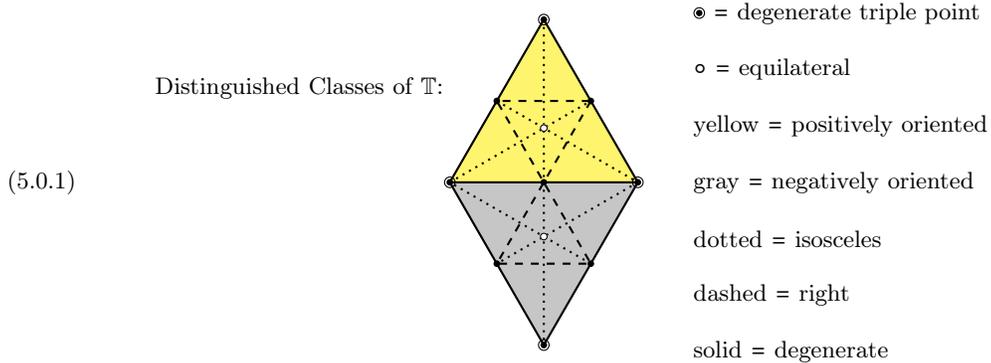
\begin{equation}\label{diamond}
\begin{minipage}{.8\textwidth}
\begin{tikzpicture}[scale=2.5]
\coordinate (A) at (-1/2,0); 
\coordinate (B) at (0,{sqrt(3)/2}); 
\coordinate (C) at (1/2,0); 
\coordinate (D) at (0,{-sqrt(3)/2}); 
\coordinate (AB) at (-1/4,{sqrt(3)/4}); 
\coordinate (BC) at (1/4,{sqrt(3)/4});
\coordinate (AC) at (0,0);
\coordinate (AD) at (-1/4,{-sqrt(3)/4});
\coordinate (CD) at (1/4,{-sqrt(3)/4});
\draw[thick] (A)--(B)--(C)--(D)--(A)--(C);
\fill (A) circle (.5pt);
\fill (B) circle (.5pt);
\fill (C) circle (.5pt);
\fill (D) circle (.5pt);
\fill (AC) circle (.5pt);
\fill (AB) circle (.5pt);
\fill (BC) circle (.5pt);
\fill (AD) circle (.5pt);
\fill (CD) circle (.5pt);
\draw[thick,dashed] (BC)--(AD)--(CD)--(AB)--(BC);
\draw[thick,dotted] (B)--(D);
\draw[thick,dotted] (CD)--(A)--(BC);
\draw[thick,dotted] (AB)--(C)--(AD);
\begin{scope}[on background layer]
\draw[fill=darkgray!30] (A)--(D)--(C)--(A);
\draw[fill=yellow!70] (A)--(B)--(C)--(A);
\end{scope}
\filldraw[white] (0,{sqrt(3)/6}) circle (.5pt);
\filldraw[white] (0,{-sqrt(3)/6}) circle (.5pt);
\draw (0,{sqrt(3)/6}) circle (.5pt);
\draw (0,{-sqrt(3)/6}) circle (.5pt);
\draw (A) circle (.8pt);
\draw (B) circle (.8pt);
\draw (C) circle (.8pt);
\draw (D) circle (.8pt);
\node[right] at (.75,.9) {$\tikz{\draw(0,0) circle (2pt);\filldraw (0,0) circle (1pt);}$ $=$ degenerate triple point};
\node[right] at (.75,.6) {$\circ$ $=$ equilateral};
\node[right] at (.75,.3) {yellow $=$ positively oriented};
\node[right] at (.75,0) {gray $=$ negatively oriented};
\node[right] at (.75,-.3) {dotted $=$ isosceles};
\node[right] at (.75,-.6) {dashed $=$ right};
\node[right] at (.75,-.9) {solid $=$ degenerate};
\node at (-1.3,.5) { Distinguished Classes of $\T$:};
\end{tikzpicture}
\end{minipage}
\end{equation}

Let $I_X$, $D_X$, $R_X$ for $X=A,B,C$ denote the points of $\T$
corresponding to 
isosceles, degenerate, and right triangles, with distinguished vertex $X$.
The standard triangle classes have the following algebraic interpretations. 
\begin{enumerate}[\rm(a)]
\item
Each $I_X\leq\mbb T$ is a subgroup isomorphic to $\R/\pi$, forming a noncontractible loop.
They are:
\[I_A=\{(-2\beta,\beta,\beta)\},\quad I_B=\{(\alpha,-2\alpha,\alpha)\},\quad
I_C=\{(\alpha,\alpha,-2\alpha)\}\]
\item
Each $ D_X\leq\mbb T$ is a subgroup isomorphic to $\R/\pi$.
They are:
\[ D_A=\{(0,\beta,-\beta)\},\quad
 D_B=\{(\alpha,0,-\alpha)\},\quad  D_C=\{(\alpha,-\alpha,0)\}\]
\item
Each $ R_X\subset\mbb T$ is a coset of $ D_X$, forming a noncontractible loop.
$ R_X$ consists of the
elements of order two $\pmod{D_X}$:
\begin{align*}
 R_A&=(\tfrac\pi 2,\tfrac\pi 4,\tfrac\pi 4)+ D_A=\{(\tfrac\pi 2,\tfrac\pi 4+\beta,\tfrac\pi 4-\beta)\}\\
 R_B&=(\tfrac\pi 4,\tfrac\pi 2,\tfrac\pi 4)+ D_B=\{(\tfrac\pi 4+\alpha,\tfrac\pi 2,\tfrac\pi 4-\alpha)\}\\
 R_C&=(\tfrac\pi 4,\tfrac\pi 4,\tfrac\pi 2)+ D_C=\{(\tfrac\pi 4+\alpha,\tfrac\pi 4-\alpha,\tfrac\pi 2)\}
\end{align*}
\item
The positive and negative equilateral classes are inverses, each of order $3$: 
$\pm(\tfrac\pi 3,\tfrac\pi 3,\tfrac\pi 3)$ in 
the intersection $I_A\cap I_B\cap I_C$.
\item
The three degenerate right-isosceles classes are of order $2$, given by
\[\{(0,\tfrac\pi 2,\tfrac\pi 2),(\tfrac\pi 2,0,\tfrac\pi 2),(\tfrac\pi 2,\tfrac\pi 2,0)\}\]
Each is the unique order-two element of the subgroups $ D_A, D_B, D_C$,
respectively.
\item
The six nondegenerate right-isosceles classes
\[\left\{\pm(\tfrac\pi 2,\tfrac\pi 4,\tfrac\pi 4),\pm(\tfrac\pi 4,\tfrac\pi 2,\tfrac\pi 4),
\pm(\tfrac\pi 4,\tfrac\pi 4,\tfrac\pi 2)\right\}\]
are generators of three cyclic groups of order $4$.
Each group contains one of the three degenerate non-equiangular isosceles class (of order $2$).
\end{enumerate}

\section{\Large Ratios}
The natural metric on 
our model computes ratios of different triangle types,
easily read off of Figure~\ref{diamond}.
Since the acute classes are inside the dashed triangles, 
the obtuse to acute ratio is $(O:A)=(3:1)$.
Normalizing the side-length of the yellow equilateral triangle to $2$, we
find the isosceles classes measure $|I|=6\sqrt 3$,
the right classes $|R|=6$,
the acute-isosceles and obtuse-isosceles $|AI|=|OI|=3\sqrt 3$,
and the degenerates $|D|=6$.
Thus we have
\[(I:R)=(2\sqrt 3:1)\,,\qquad(I:OI)=(I:AI)=(2:1)\,,\qquad(D:R)=(1:1).\]
We do not think this model is any more ``metrically correct'' than other models, such
as the sphere in \cite{Beh} and \cite{ES15}. In fact, we believe there are strong arguments that both are incorrect,
essentially because both omit important degenerate classes,
as mentioned in Subsection~\ref{missing}.

\section{\Large Group Action}\label{groupoid}

A nondegenerate scalene triangle has a single absolute similarity class,
but twelve labeled, oriented similarity classes:
six corresponding to the ways of
assigning the three angles assigned to the three vertices,
and two for each orientation.
In general $\T$ admits a group action by $D_6=\br{r,s}$ whose orbits are classes that
are {\it absolutely similar}, i.e., similar as unoriented, unlabeled 
inscribable classes.
This is depicted in Figure~\ref{hex}:

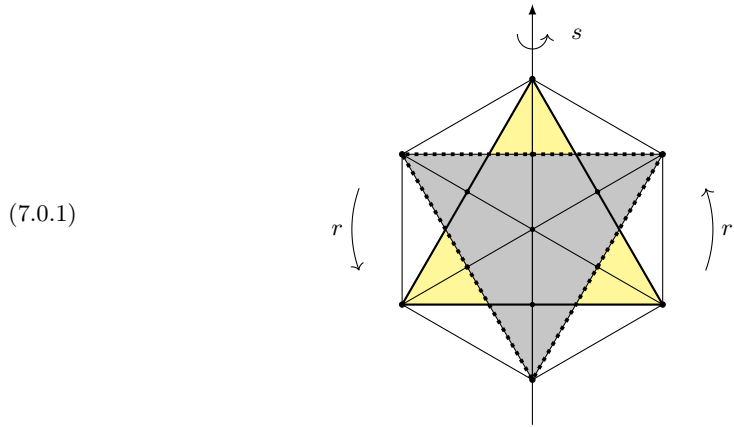
\begin{equation}\label{hex}
\begin{minipage}{.4\textwidth}
\begin{tikzpicture}[scale=2]
\coordinate (O) at (0,0); \fill (O) circle (.5pt);
\coordinate (A) at ({cos(30)},{sin(30)}); \filldraw (A) circle (.5pt);
\coordinate (B) at ({cos(90)},{sin(90)}); \filldraw (B) circle (.5pt);
\coordinate (C) at ({cos(150)},{sin(150)}); \filldraw (C) circle (.5pt);
\coordinate (D) at ({cos(210)},{sin(210)}); \filldraw (D) circle (.5pt);
\coordinate (E) at ({cos(270)},{sin(270)}); \filldraw (E) circle (.5pt);
\coordinate (F) at ({cos(330)},{sin(330)}); \filldraw (F) circle (.5pt);
\coordinate (AC) at ({(1/2)*cos(30)},{(1/2)*sin(30)}); \fill (AC) circle (.5pt);
\coordinate (CE) at ({(1/2)*cos(90)},{(1/2)*sin(90)}); \fill (CE) circle (.5pt);
\coordinate (EA) at ({(1/2)*cos(150)},{(1/2)*sin(150)}); \fill (EA) circle (.5pt);
\coordinate (BD) at ({(1/2)*cos(210)},{(1/2)*sin(210)}); \fill (BD) circle (.5pt);
\coordinate (DF) at ({(1/2)*cos(270)},{(1/2)*sin(270)}); \fill (DF) circle (.5pt);
\coordinate (FB) at ({(1/2)*cos(330)},{(1/2)*sin(330)}); \fill (FB) circle (.5pt);
\draw[thick] (B)--(D)--(F)--(B); 
\draw[line width=0.5mm, dotted] (A)--(C)--(E)--(A); 
\draw[thin] (A)--(B)--(C)--(D)--(E)--(F)--(A); 
\draw[thin] (A)--(D);
\draw[thin] (B)--(E);
\draw[thin] (C)--(F);
\draw[thin,->,-latex] (0,-1.3)--(0,1.5); 
\draw[domain=-180:0, ->] plot ({.1*cos(\x)},{1.3+.1*sin(\x)});
\draw[domain=-20:20, ->] plot ({.4+.8*cos(\x)},{.8*sin(\x)});
\draw[domain=160:200, ->] plot ({-.4+.8*cos(\x)},{.8*sin(\x)});
\node[right] at (.2,1.3) {$s$};
\node[right] at (1.2,0) {$r$};
\node[left] at (-1.2,0) {$r$};
\begin{scope}[on background layer]
\draw[fill=yellow!50] (B)--(D)--(F)--(B);
\draw[fill=darkgray!30] (A)--(C)--(E)--(A);
\end{scope}
\end{tikzpicture}
\end{minipage}
\end{equation}
The explicit action on a general point is given by $r(\alpha,\beta,\gamma)=(-\beta,-\gamma,-\alpha)$
and $s(\alpha,\beta,\gamma)=(\beta,\alpha,\gamma)$.
The resulting transformation groupoid is $\T\times D_6$.
The bipyramid of \cite[Exercise 1.38]{Beh} has the same groupoid, and the construction is parallel.
The quotient $\T/D_6$ is the (coarse) moduli set of absolute possibly-degenerate inscribable classes,
drawn in Figure~\ref{stacky} with ``stacky'' multiplicities, which are the orders of the stabilizer subgroups.
\begin{equation}\label{stacky}
\begin{minipage}{.9\textwidth}
\begin{tikzpicture}
\coordinate (A) at (0,0);
\coordinate (B) at (0,1);
\coordinate (C) at ({-sqrt(3)},0);
\coordinate (D) at ({-sqrt(3)/3},2/3);
\draw[thick] (A)--(B)--(C)--(A);
\filldraw[white] (B) circle (1.5pt);
\fill (A) circle (1.5pt);
\draw (B) circle (1.5pt) node[above right] {equilateral $(\times 6)$};
\fill (C) circle (1.5pt);
\draw (-14/12,1/3) arc (0:90:1 and .8) node[left] {isosceles $(\times 2)$};
\draw (0,2/3) arc (0:73:3 and .5); 
\draw (-14/12,1/3) circle (1.5pt);
\draw (0,2/3) circle (1.5pt);
\draw (-7/8,0) circle (1.5pt);
\draw (-7/8,0) arc (180:73:.5 and -.4) node[right] {double point $(\times 2)$};
\draw (A)--(.4,.2) node[right] {right double point\; $(\times 4)$}; 
\draw (C)--(-1.75-.4,.2) node[left] {doubled simple point\; $(\times 12)$};
\begin{scope}[on background layer]
\draw[pattern=north west lines] (A)--(B)--(C)--(A);
\end{scope}
\node at (-7/8,-1.1) { $\T/D_6$ with Stacky Multiplicities};
\end{tikzpicture}
\end{minipage}
\end{equation}
A doubled simple point is a double point with interior angles all
zero $\pmod\pi$, and a right double point is a double point with two right angles.
To say $\T/D_6$ is a moduli set simply is to say there is a bijection between its points and
the (absolute) inscribable classes.

\section{\Large Moduli Space and Quotient Stack}\label{moduliquotient}

We will now view $\T$ as a variety as per Corollary~\ref{algvar},
and prove it is a proper moduli space of labeled, oriented, possibly-degenerate inscribable classes
in the algebrogeometric setting.
We use the term {\it $\R$-variety} to mean a nonsingular, integral, separated scheme of finite type over $\R$.
Every $\R$-variety is a smooth manifold over $\R$. We write $\Var{\R}$ for the category of
$\R$-varieties, and in the following assume some elementary background in algebraic geometry
and category theory.

\begin{Definition}\label{families}\rm
If $B$ is an $\R$-variety, let $\pi_\T$ and $\pi_B$ be the projections on $\T\times_\R B$.
\begin{enumerate}[\rm (a)]
\item\label{a}
A {\it family in $\T$ over $B$} is a subvariety $F\subset \T\times_\R B$ such that
$\pi_B|_F:F\to B$ is an isomorphism.
\item\label{b}
The {\it pullback} $\phi^*(F/B)$ along a morphism $\phi:C\to B$
is the family $F_C:=(F\times_B C)/C$.
\item
The (contravariant) \textit{moduli functor}
$\FF_{\,\T}:\Var{\R}\to\Set$
assigns to each $B\in\Var{\R}$ the set of families $F/B$ in $\T$, and
to each map $\phi:C\to B$ in $\Var{\R}$ the pullback of \eqref{b}.
\end{enumerate}
\end{Definition}

\begin{Remark}\label{alt}
Alternatively, a family is given by a morphism $f=(\alpha,\beta,\gamma):B\lr (\R/\pi)^3$ in $\Var{\R}$
whose image lies on $\T$.
That is, we assign to each $b\in B$ an
ordered triple $f(b)=(\alpha(b),\beta(b),\gamma(b))$, where $\alpha,\beta,\gamma:B\to\R/\pi$
are (regular) morphisms.
This phrasing
emphasizes the requirement that families vary continuously (regularly) over the points of $B$,
and is more in line with the definitions of families given in \cite{Beh} and other sources.
But it's the same:
Since $\T/\R$ is separated, the graph $F=\Gamma_f$ is closed in $\T\times_\R B$, and
since $f$ is regular, $\pi_B|_F:F\to B$ is an isomorphism.
Conversely,
a family $F/B$ in \eqref{a}
defines a morphism $f:B\to\T$ by composing $(\pi_B|_F)^{-1}:B\to F$ 
with $\pi_\T:\T\times B\to\T$. This is a morphism since $\pi_B|_F$ is an isomorphism.
\end{Remark}

\begin{Theorem}\label{Theta}
The functor $\FF_{\,\T}$ is represented by $\T$. That is, it is isomorphic to $\Hom_{\Var{\R}}(-,\T)$.
\end{Theorem}

\begin{proof}
More generally if a set of objects is 
parameterized by an $S$-scheme $X$, then the functor that assigns
to each $S$-scheme $B$ the set of families of objects over $B$,
is representable by $X$ (see \cite{Beh}). 

We construct a natural transformation $\Theta:\FF_{\,\T}\to\Hom_{\Var{\R}}(-,\T)$.
A family $F\in\FF_{\,\T}(B)$ defines a morphism $f:B\to\T$ by Remark~\ref{alt}.
Define the component of $\Theta$ on the (arbitrary) base $B$ by $\Theta_B(F)=f$.
To prove $\Theta$ is a natural transformation we must show that a morphism
$\phi:C\to B$ determines a commutative diagram
\[\begin{tikzcd}[ampersand replacement=\&]
{\FF_{\,\T}(B)} \&\& {\Hom_{\Var{\R}}(B,\T)} \\
\\
{\FF_{\,\T}(C)} \&\& {\Hom_{\Var{\R}}(C,\T)}
\arrow["{\Theta_B}", from=1-1, to=1-3]
\arrow["\FF_{\,\T} (\phi)"', from=1-1, to=3-1]
\arrow["{\Hom_{\Var{\R}}(\phi,\T)}", from=1-3, to=3-3]
\arrow["{\Theta_C}"', from=3-1, to=3-3]
\end{tikzcd}\]
Suppose $F\in\FF_{\,\T}(B)$ defines the $B$-point $f:B\to\T$, as above.
Then $\Theta_B(F)=f$, and by definition $(\Hom_{\Var{\R}}(\phi,\T)\circ\Theta_B)(F)=f\circ\phi$.
On the other hand,
by definition $\FF_{\,\T}(\phi)(F)=F\times_B C=F_C$,
and $\Theta_C(F_C)$ is the composition $\pi_\T\circ(\pi_C|_{F_C})^{-1}$.
Suppose $(\eta,c)$ is a point of $F_C$. By definition of fiber product
$\pi_B(\eta)=\phi(c)$, so $(\eta,c)=((\pi_\T(\eta),\phi(c)),c)$ is the corresponding
point on $\T\times_\R B\times_B C$. 
Therefore $\Theta_C(F_C)(c)\df\pi_\T\circ(\pi_C|_{F_C})^{-1}(c)=\pi_\T(\eta)$,
and $\pi_\T(\eta)=f(\phi(c))$ by definition of $f$.
We conclude $\Theta_C\circ\FF_{\,\T}(\phi)=f\circ\phi$.
Since $F/B$ was arbitrary,
this shows the diagram commutes, hence
$\Theta$ is a natural transformation.

We invert $\Theta$ by taking $f\in\Hom_{\Var{\R}}(B,\T)$ to
the family $F/B\in\FF_{\,\T}(B)$ defined in Remark~\ref{alt}.
The proof that this inverts $\Theta$ is straightforward, and 
shows $\FF_{\,\T}$ and $\Hom_{\Var{\R}}(-,\T)$ are isomorphic.
\end{proof}

\begin{Corollary}\label{moduli}
$\T$ is a fine moduli space for the set of all labeled, oriented,
possibly-degenerate inscribable triangle classes,
which is a compactification of $\T^\circ=\T-\{\alpha\beta\gamma=0\}$, 
the fine moduli space of labeled, oriented, nondegenerate inscribable classes.
In particular, the family $[\Delta_{\rm insc}]=\T$ is a universal family, given explicitly by $\sigma$
in Theorem~\ref{map}.
\end{Corollary}

\begin{proof}
Already there is a bijection between labeled, oriented, possibly-degenerate inscribable classes
and points of $\T$, and it remains to verify $\Hom_{\Var{\R}}(-,\T)$ has the universal property with
respect to natural transformations $\FF_{\,\T}\to\Hom_{\Var{\R}}(-,Y)$ for $Y\in\Var{\R}$.
But $\FF_{\,\T}$ and $\Hom_{\Var{\R}}(-,\T)$ are isomorphic via $\Theta$ of Theorem~\ref{Theta}, 
so for any ${\Psi}:\FF_{\,\T}$ $\to\Hom_{\Var{\R}}(-,Y)$
there is trivially a natural transformation $\Phi:\Hom_{\Var{\R}}(-,\T)\to\Hom_{\Var{\R}}(-,Y)$
such that $\Psi=\Phi\circ\Theta$. This proves $\T$ is a fine moduli space.
That it compactifies $\T^\circ$ is trivial since it is compact and it is the closure,
and that $\T^\circ$ is a fine moduli space is also immediate via the open immersion $\T^\circ\subset\T$.
The diagonal family of $D\subset\T\times_\R\T$ in $\FF_{\,\T}(\T)$ maps to $\id_\T\in\Hom_{\Var{\R}}(\T,\T)$,
and is the accompanying universal family.
\end{proof}

\subsubsection{\bf Moduli Stacks}
The coarse moduli set $\T/D_6$ of absolute (unlabeled, 
unoriented) possibly-degenerate inscribable classes in Figure~\ref{stacky} is not an $\R$-variety,
hence not a fine moduli space, 
essentially because the symmetry group $D_6$ does not act freely on $\T$, as indicated by the nontrivial
stacky multiplicities labeled in \eqref{stacky}.
More intrinsically, the stacky points cause problems because the corresponding
nontrivial automorphisms of individual classes can be
used to construct non-isomorphic families that the moduli set $\T/D_6$ doesn't distinguish
(see \cite{Beh}), which thwarts the construction of a universal family of absolute classes.

Nevertheless, 
since the moduli set of absolute inscribable classes is the orbit space of
the $\R$-variety $\T$ under $D_6$, it has the
geometric structure of a quotient stack, denoted $[\T/D_6]$
(\cite[Example 4.8]{DM69}, \cite[1.24]{Beh}, \cite[0.6.5]{Alper}).
By definition, $[\T/D_6]$ 
is the category whose objects are $D_6$-torsors over $\R$-varieties $B$ 
that are equipped with $D_6$-equivariant maps into $\T$.
These are given by diagrams
\[\begin{tikzcd}[ampersand replacement=\&]
B'\&\T\\ B
\arrow[r,"f",from=1-1,to=1-2,-latex]
\arrow[d,from=1-1,to=2-1,-latex]
\end{tikzcd}\]
denoted $(B'/B,f)$,
where $B'/B$ is a $D_6$-torsor and $f$ is $D_6$-equivariant.
The morphisms $(\phi,\theta):(C'/C,g)\to(B'/B,f)$ are pairs $(\phi,\theta)$ of morphisms $\phi:C\to B$ in $\Var{\R}$
and $D_6$-equivariant isomorphisms $\theta:C'\to B'\times_B C\to B'$ satisfying $g=f\circ\theta$:
\[\begin{tikzcd}[ampersand replacement=\&]
C'\&B'\&\T\\ C \& B
\arrow[r,"\theta",from=1-1,to=1-2,-latex]
\arrow[r,"f",from=1-2,to=1-3,-latex]
\arrow[r,"g",bend left,from=1-1,to=1-3,-latex]
\arrow[d,from=1-1,to=2-1,-latex]
\arrow[d,from=1-2,to=2-2,-latex]
\arrow[r,"\phi",from=2-1,to=2-2,-latex]
\end{tikzcd}\]
By incorporating the automorphisms of individual classes into each family over $B$,
the stack distinguishes nonisomorphic families. Though it is not a variety, $[\T/D_6]$ is as close
to a fine moduli space as we get.

\bibliographystyle{alpha} 
\bibliography{../../MathDocs/mathdocs.bib}

\end{document}